\documentclass[12pt]{article}
\usepackage{amsmath}
\usepackage{graphicx,psfrag,epsf}
\usepackage{enumerate}
\usepackage{natbib}
\usepackage{url} 

\usepackage{bm}
\usepackage{amsfonts}
\usepackage{subcaption}
\usepackage{multirow}

\newcommand{\blind}{0}

\newtheorem{theorem}{Theorem}

\newtheorem{proposition}{Proposition}

\newtheorem{remark}{Remark}

\newcommand{\eop}{\hfill $\Box$ \\ \\}
\newcommand{\Cov}{\mathrm{Cov}}

\DeclareMathOperator{\sign}{sign}

\begin{document}

\def\spacingset#1{\renewcommand{\baselinestretch}%
{#1}\small\normalsize} \spacingset{1}


\if0\blind
{
  \title{\bf Fused Lasso Nearly Isotonic Signal Approximation in General Dimensions}
  \author{Vladimir Pastukhov\\
    Department of Computer Science and Engineering,\\ Chalmers, Sweden\\
    Email: vlapas@chalmers.se
}
  \maketitle
} \fi

\if1\blind
{
  \bigskip
  \bigskip
  \bigskip
  \begin{center}
    {\LARGE\bf Title}
\end{center}
  \medskip
} \fi

\bigskip
\begin{abstract}
In this paper we introduce and study fused lasso nearly-isotonic signal approximation, which is a combination of fused lasso and generalized nearly-isotonic regression. We show how these three estimators relate to each other, derive solution to the general problem, show that it is computationally feasible and provides a trade-off between piecewise monotonicity, sparsity and goodness-of-fit. Also, we derive an unbiased estimator of the degrees of freedom of the approximator. 
\end{abstract}

\noindent%
{\it Keywords:}  Constrained inference; Isotonic regression; Nearly-isotonic Regression; Fused lasso
\vfill

\newpage
\spacingset{1.45} 
\section{Introduction}
\label{sec:intro}

This work is motivated by recent papers in nearly-constrained estimation and papers in penalized least squared regression. The subject of penalized estimators starts with $L_{1}$-penalisation  \citep{tibshirani1996regression}, i.e.
\begin{equation*}\label{*}
\hat{\bm{\beta}}^{L}(\bm{y}, \lambda_{L}) = \underset{\bm{\beta} \in \mathbb{R}^{n}}{\arg \min} \, \frac{1}{2} ||\bm{y} - \bm{\beta}||_{2}^{2} + \lambda_{L} ||\bm{\beta}||_{1},
\end{equation*}
which is called lasso signal approximation. The $L_{2}$-penalisation is usually addressed as ridge regression \citep{hoerl1970ridge} or sometimes as Tikhonov-Philips regularization \citep{phillips1962technique, tikhonov1995numerical}.

First, to set up the ideas and for simplicity of notation we consider one dimensional cases of the penalized estimators. In the next subsection we generalise these estimators to the case of isotonic constraints with respect to a general partial order. 

For a given sequence of data points $\bm{y} \in \mathbb{R}^{n}$ the fusion approximator (cf. \cite{rinaldo2009properties}) is given by
\begin{equation}\label{F}
\hat{\bm{\beta}}^{F}(\bm{y}, \lambda_{F}) = \underset{\bm{\beta} \in \mathbb{R}^{n}}{\arg \min} \, \frac{1}{2} ||\bm{y} - \bm{\beta}||_{2}^{2} +\lambda_{F} \sum_{i=1}^{n-1}|\beta_{i} - \beta_{i+1}|.
\end{equation}
Further, the combination of fusion estimator and lasso is called fused lasso estimator and is given by:
\begin{equation}\label{FL}
\hat{\bm{\beta}}^{FL}(\bm{y}, \lambda_{F}, \lambda_{L}) = \underset{\bm{\beta} \in \mathbb{R}^{n}}{\arg \min} \, \frac{1}{2} ||\bm{y} - \bm{\beta}||_{2}^{2} +\lambda_{F}\sum_{i=1}^{n-1}|\beta_{i} - \beta_{i+1}| + \lambda_{L} ||\bm{\beta}||_{1}.
\end{equation}
The fused lasso was introduced in cf. \citet{tibshirani2005sparsity} and its asymptotic properties were studied in detail in \citet{rinaldo2009properties}.
\begin{remark}\label{}
In the paper \citet{tibshirani2011solution} the estimator in (\ref{F}) is called  the fused lasso, while the estimator in (\ref{FL}) is addressed as the sparse fused lasso.
\end{remark}

In the area of constrained inference the basic and simplest problem is the isotonic regression in one dimension. For a given sequence of data points $\bm{y} \in \mathbb{R}^{n}$ the isotonic regression is the following approximation
\begin{equation}\label{orMLEp}
\hat{\bm{\beta}}^{I} = \underset{\bm{\beta} \in \mathbb{R}^{n}}{\arg\min}||\bm{y} - \bm{\beta}||_{2}^{2}, \quad \text{subject to} \quad \beta_{1} \leq \beta_{2} \leq \dots \leq \beta_{n}, 
\end{equation}
i.e. the isotonic regression is the $\ell^{2}$-projection of the vector $\bm{y}$ onto the set of non-increasing vectors in $\mathbb{R}^{n}$. The notion of isotonic "regression" in this context might be confusing. Nevertheless, it is a standard notion in this subject, cf., for example, the papers \citet{best1990active, stout2013isotonic}, where the notation "isotonic regression" is used for the isotonic projection of a general vector. Also, in this paper we use notations "regression", "estimator" and "approximator" interchangeably. 

The nearly-isotonic regression, introduced in \cite{tibshirani2011nearly} and studied in detail in \cite{minami2020estimating}, is a less restrictive version of the isotonic regression and is given by the following optimization problem
\begin{equation}\label{NI}
\hat{\bm{\beta}}^{NI}(\bm{y}, \lambda_{NI}) = \underset{\bm{\beta} \in \mathbb{R}^{n}}{\arg \min} \, \frac{1}{2} ||\bm{y} - \bm{\beta}||_{2}^{2} + \lambda_{NI}\sum_{i=1}^{n-1}|\beta_{i} - \beta_{i+1}|_{+},
\end{equation}
where $x_{+} =  x \cdot 1\{x > 0 \}$.

In this paper we combine fused lasso estimator with nearly-isotonic regression and call the resulting estimator as \textit{fused lasso nearly-isotonic signal approximator}, i.e. for a given sequence of data points $\bm{y} \in \mathbb{R}^{n}$ the problem in one dimensional case is the following optimization
\begin{equation}\label{FLNI}
\begin{aligned}
\hat{\bm{\beta}}^{FLNI}(\bm{y}, \lambda_{F}, \lambda_{L}, \lambda_{NI}) = {} & \\
 \underset{\bm{\beta} \in \mathbb{R}^{n}}{\arg \min} \, \frac{1}{2} ||\bm{y} - \bm{\beta}||_{2}^{2} +& \lambda_{F} \sum_{i=1}^{n-1}|\beta_{i} - \beta_{i+1}| + \lambda_{L} ||\bm{\beta}||_{1} + \lambda_{NI}\sum_{i=1}^{n-1}|\beta_{i} - \beta_{i+1}|_{+}.
\end{aligned}
\end{equation}
Also, in the case of $\lambda_{F} \neq 0$ and $\lambda_{NI} \neq 0$, with $\lambda_{L} = 0$,  we call the estimator as \textit{fused nearly-isotonic regression}, i.e.
\begin{equation}\label{FNI}
\begin{aligned}
\hat{\bm{\beta}}^{FNI}(\bm{y}, \lambda_{F}, \lambda_{NI}) \equiv {} & \hat{\bm{\beta}}^{FLNI}(\bm{y}, \lambda_{F}, 0, \lambda_{NI}) = \\ &\underset{\bm{\beta} \in \mathbb{R}^{n}}{\arg \min} \, \frac{1}{2} ||\bm{y} - \bm{\beta}||_{2}^{2}  + \lambda_{F} \sum_{i=1}^{n-1}|\beta_{i} - \beta_{i+1}| + \lambda_{NI}\sum_{i=1}^{n-1}|\beta_{i} - \beta_{i+1}|_{+}.
\end{aligned}
\end{equation}
The generalisation of nearly-isotonic regression in (\ref{FNI}) was proposed in the conclusion of the paper \citet{tibshirani2011nearly}. 

We will state the problem in (\ref{FLNI}) for the case of isotonic constraints with respect to a general partial order, but, first, we have to introduce the notation. 

\subsection{Notation}
We start with basic definitions of partial order and isotonic regression.  Let $\mathcal{I} = \{\bm{i}_{1}, \dots, \bm{i}_{n}\}$ be some index set. Next, we define the following binary relation $\preceq$ on $\mathcal{I}$.

A binary relation $\preceq$ on $\mathcal{I}$ is called partial order if 
\begin{itemize}
\item it is reflexive, i.e. $\bm{j}\preceq\bm{j}$ for all $\bm{j} \in \mathcal{I}$;
\item it is transitive, i.e. $\bm{j}_{1}, \bm{j}_{2}, \bm{j}_{3} \in \mathcal{I}$, $\bm{j}_{1} \preceq \bm{j}_{2}$ and $\bm{j}_{2} \preceq \bm{j}_{3}$ imply $\bm{j}_{1} \preceq \bm{j}_{3}$;
\item it is antisymmetric, i.e. $\bm{j}_{1}, \bm{j}_{2} \in \mathcal{I}$, $\bm{j}_{1} \preceq \bm{j}_{2}$ and $\bm{j}_{2} \preceq \bm{j}_{1}$ imply $\bm{j}_{1} = \bm{j}_{2}$.
\end{itemize}

Further, a vector $\bm{\beta}\in\mathbb{R}^{n}$ indexed by $\mathcal{I}$ is called isotonic with respect to the partial order $\preceq$ on $\mathcal{I}$ if $\bm{j}_{1} \preceq \bm{j}_{2}$ implies $\beta_{\bm{j}_{1}} \leq \beta_{\bm{j}_{2}}$. We denote the set of all isotonic vectors in $\mathbb{R}^{n}$ with respect to the partial order $\preceq$ on $\mathcal{I}$ by $\bm{\mathcal{B}}^{is}$, which is also called isotonic cone. Next, a vector $\bm{\beta}^{I}\in \mathbb{R}^{n}$ is isotonic regression of an arbitrary vector $\bm{y} \in \mathbb{R}^{n}$ over the pre-ordered index set $\mathcal{I}$ if 
\begin{eqnarray}\label{Ipo}
\bm{\beta}^{I} = \underset{\bm{\beta} \in \bm{\mathcal{B}}^{is}}{\arg \min} \sum_{\bm{j} \in \mathcal{I}}(\beta_{\bm{j}} - y_{\bm{j}})^{2}.
\end{eqnarray}

For any partial order relation $\preceq$ on $\mathcal{I}$ there exists directed graph $G = (V,E)$, with $V = \mathcal{I}$ and 
\begin{eqnarray}
E = \{(\bm{j}_{1}, \bm{j}_{2}), \, \text{where} \, \bm{j_{1}} \, \text{and} \, \bm{j_{2}} \,\text{are certain elements of} \, \mathcal{I}\},
\end{eqnarray}
such that an arbitrary vector $\bm{\beta} \in \mathbb{R}^{n}$ is isotonic with respect to $\preceq$ iff $\beta_{\bm{j_{1}}} \leq \beta_{\bm{j_{2}}}$ for any $(\bm{j}_{1}, \bm{j}_{2}) \in E$. Therefore, equivalently to the definition in (\ref{Ipo}), a vector $\bm{\beta}^{I}\in \mathbb{R}^{n}$ is isotonic regression of an arbitrary vector $\bm{y} \in \mathbb{R}^{n}$  indexed by the partially ordered index set $\mathcal{I}$ if
\begin{eqnarray}\label{IpoG}
\bm{\beta}^{I} = \underset{\bm{\beta}}{\arg \min} \sum_{\bm{j} \in \mathcal{I}}(\beta_{\bm{j}} - y_{\bm{j}})^{2}, \quad \text{subject to} \quad \beta_{\bm{j_{1}}} \leq \beta_{\bm{j_{2}}}, \quad \text{whenever} \quad (\bm{j_{1}}, \bm{j_{2}}) \in E.
\end{eqnarray}

Further, for the directed graph $G = (V, E)$, which corresponds to the partial order $\preceq$ on $\mathcal{I}$, the nearly-isotonic regresson of $\bm{y}\in \mathbb{R}^{n}$ indexed by $\mathcal{I}$ is given by
\begin{equation}\label{NIG}
\hat{\bm{\beta}}^{NI}(\bm{y}, \lambda_{NI}) = \underset{\bm{\beta} \in \mathbb{R}^{n}}{\arg \min} \, \frac{1}{2} ||\bm{y} - \bm{\beta}||_{2}^{2} + \lambda_{NI}\sum_{(\bm{i},\bm{j})\in E}|\beta_{\bm{i}} - \beta_{\bm{j}}|_{+}.
\end{equation}
This generalisation of nearly-isotonic regression was introduced and studied in \citet{minami2020estimating}. 

Next, fussed and fussed lasso approximators for a general directed graph $G = (V, E)$ are given by
\begin{equation}\label{FG}
\hat{\bm{\beta}}^{F}(\bm{y}, \lambda_{F}) = \underset{\bm{\beta} \in \mathbb{R}^{n}}{\arg \min} \, \frac{1}{2} ||\bm{y} - \bm{\beta}||_{2}^{2} +\lambda_{F}\sum_{(\bm{i},\bm{j})\in E}|\beta_{\bm{i}} - \beta_{\bm{j}}|,
\end{equation}
and
\begin{equation}\label{FLG}
\hat{\bm{\beta}}^{FL}(\bm{y}, \lambda_{F}, \lambda_{L}) = \underset{\bm{\beta} \in \mathbb{R}^{n}}{\arg \min} \, \frac{1}{2} ||\bm{y} - \bm{\beta}||_{2}^{2} +\lambda_{F}\sum_{(\bm{i},\bm{j})\in E}|\beta_{\bm{i}} - \beta_{\bm{j}}| + \lambda_{L} ||\bm{\beta}||_{1}.
\end{equation}
These optimization problems were introduced and solved for a general directed graph in \citet{tibshirani2011solution}.

Further, let $D$ denote the oriented incidence matrix for the directed graph $G = (V, E)$, corresponding to $\preceq$ on $\mathcal{I}$ . We choose the orientation of $D$ in the following way. Assume that the graph $G$ with $n$ vertexes has $m$ edges. Next, assume we label the vertexes by $\{1, \dots, n\}$ and edges by $\{1, \dots, m\}$. Then $D$ is $m\times n$ matrix with
\begin{equation*}\label{}
D_{i,j} =   \begin{cases}
   1, & \quad \text{if vertex $j$ is the source of edge $i$} , \\
    -1, & \quad \text{if vertex $j$ is the target of edge $i$},\\
    0, & \quad \text{otherwise}.
  \end{cases}
\end{equation*}

\begin{figure}[h!]
\centering
\begin{minipage}[b]{.45\textwidth}
\centering
\includegraphics[width=\textwidth]{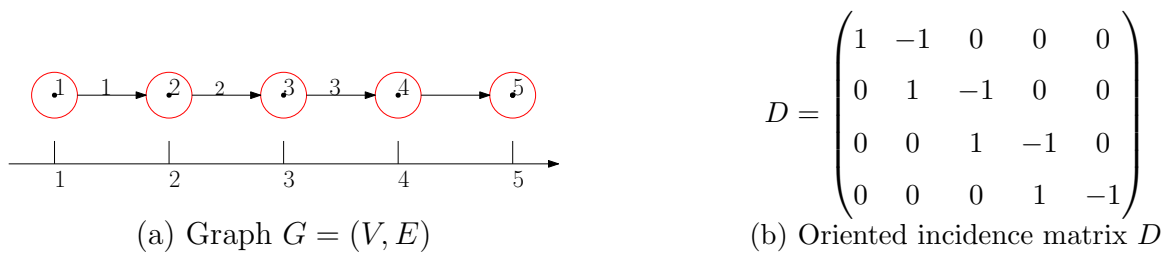}\\
(a) Graph $G=(V,E)$
\end{minipage}\hfill
\begin{minipage}[b]{.45\textwidth}
\centering
\small
\(
D = \begin{pmatrix}
1 & -1 & 0 & 0 & 0  \\
0 & 1 & -1 & 0 & 0 \\
0 & 0 & 1 & -1 & 0 \\
0 & 0 & 0 & 1 & -1
\end{pmatrix}
\)\\
(b) Oriented incidence matrix $D$
\end{minipage}
\caption{Graph for monotonic contstraints and oriented incidence matrix}\label{mgrmtr}
\end{figure}

In order to clarify the notations we consider the following examples of partial order relations. First, let us consider the monotonic order relation in one dimensional case. Let $\mathcal{I} = \{1, \dots, n\}$, and for $j_{1} \in \mathcal{I}$ and $j_{2} \in \mathcal{I}$  we naturally define $j_{1}\preceq j_{2}$ if $j_{1} \leq j_{2}$. Further, if we let $V = \mathcal{I}$ and $E = \{(i, i+1): i = 1, \dots, n-1 \}$, then $G = (V, E)$ is the directed graph which correspond to the one dimensional order relation on $\mathcal{I}$. Figure \ref{mgrmtr} displays the graph and the incidence matrix for the graph.

\begin{figure}[h!]
\centering
\begin{minipage}[b]{.45\textwidth}
\centering
\includegraphics[width=\textwidth]{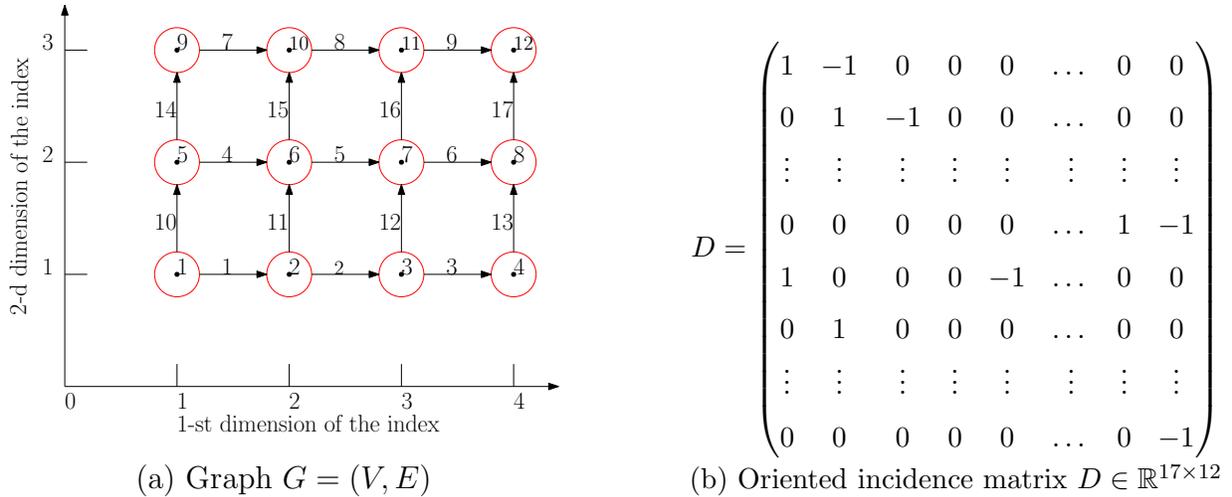}\\
(a) Graph $G=(V,E)$
\end{minipage}\hfill
\begin{minipage}[b]{.45\textwidth}
\centering
\small
\(
D = \begin{pmatrix}
1 & -1 & 0 & 0 & 0 & \dots & 0 & 0 \\
0 & 1 & -1 & 0 & 0 & \dots & 0 & 0 \\
\vdots & \vdots & \vdots & \vdots  & \vdots  & \vdots & \vdots & \vdots \\
0 & 0 & 0 & 0 & 0 & \dots & 1 & -1\\
1 & 0 & 0 &0 & -1 & \dots & 0 & 0 \\
0 & 1 & 0 &0 & 0 & \dots & 0 & 0 \\
\vdots & \vdots & \vdots  & \vdots & \vdots  & \vdots & \vdots & \vdots \\
0 & 0 & 0 & 0 & 0 & \dots & 0 & -1 
\end{pmatrix}
\)\\
(b) Oriented incidence matrix $D \in \mathbb{R}^{17\times 12}$
\end{minipage}
\caption{Graph for bimonotonic contstraints and oriented incidence matrix}\label{bimgrmtr}
\end{figure}

Next, we consider two dimensional case with bimonotonic constraints. The notion of bimonotonicity was first introduced in \citet{beran2010least} and it means the following. Let us consider the index set
\begin{eqnarray*}
\mathcal{I} = \{ \bm{i}= (i^{(1)} ,i^{(2)}): \, i^{(1)}=1,2,\dots, n_{1}, \,  i^{(2)}=1,2,\dots, n_{2}\}
\end{eqnarray*}
with the following order relation $\preceq$ on it: for $\bm{j}_{1}, \bm{j}_{2}\in \mathcal{I}$ we have $\bm{j}_{1} \preceq \bm{j}_{2}$ iff $j^{(1)}_{1} \leq j^{(1)}_{2}$ and $j^{(2)}_{1} \leq j^{(2)}_{2}$. Then, a vector $\bm{\beta}\in\mathbb{R}^{n}$, with $n=n_{1}n_{2}$, indexed by $\mathcal{I}$ is called bimonotone if it is isotonic with respect to bimonotone order $\preceq$ defined on its index $\mathcal{I}$. Further, we define the directed graph $G = (V, E)$ with vertexes $V = \mathcal{I}$, and the edges
\begin{eqnarray*}
\begin{aligned}
E ={}&  \{((l, k),(l, k+1) ): \, 1 \leq l \leq n_{1}, 1 \leq k \leq n_{2} - 1\}\\
\cup \,& \{((l, k),(l+1, k) ): \, 1 \leq l \leq n_{1}-1, 1 \leq k \leq n_{2} \}.
\end{aligned}
\end{eqnarray*}
The labeled graph and its incidence matrix are displayed on Figure \ref{bimgrmtr}.

\subsection{General statement of the problem}
Now we can state the general problem studied in this paper. Let $\bm{y} \in \mathbb{R}^{n}$ be a signal indexed by the index set $\mathcal{I}$ with the partial order relation $\preceq$ defined on $\mathcal{I}$. Next, let $G=(V,E)$ be the directed gpraph corresponding to $\preceq$ on $\mathcal{I}$. The fussed lasso nearly isotonic signal approximation with respect to $\preceq$ on $\mathcal{I}$ (or, equivalently, to the directed graph $G = (V, E)$, corresponding to $\preceq$) is given by
\begin{equation}\label{FLNIG}
\begin{aligned}
\hat{\bm{\beta}}^{FLNI}(\bm{y}, \lambda_{F}, \lambda_{L}, \lambda_{NI}) = {} & \\
 \underset{\bm{\beta} \in \mathbb{R}^{n}}{\arg \min} \, \frac{1}{2} ||\bm{y} - \bm{\beta}||_{2}^{2} +& \lambda_{F} \sum_{(\bm{i},\bm{j})\in E}|\beta_{\bm{i}} - \beta_{\bm{j}}| + \lambda_{L} ||\bm{\beta}||_{1} + \lambda_{NI}\sum_{(\bm{i},\bm{j})\in E}|\beta_{\bm{i}} - \beta_{\bm{j}}|_{+}.
\end{aligned}
\end{equation}
Therefore, the estimator in (\ref{FLNIG}) is a combination of the estimators in (\ref{NIG}) and (\ref{FLG}).

Equivalently, we can rewrite the problem in the following way:
\begin{equation}\label{FLNIGim}
\hat{\bm{\beta}}^{FLNI}(\bm{y}, \lambda_{F}, \lambda_{L}, \lambda_{NI}) = {} 
 \underset{\bm{\beta} \in \mathbb{R}^{n}}{\arg \min} \, \frac{1}{2} ||\bm{y} - \bm{\beta}||_{2}^{2} + \lambda_{F} ||D\bm{\beta}||_{1} + \lambda_{L} ||\bm{\beta}||_{1} + \lambda_{NI}||D\bm{\beta}||_{+},
\end{equation}
where $D$ is the oriented incidence matrix of the graph $G =(V, E)$.

Analogously to the definition in one dimensional case, if $\lambda_{L} = 0$ we call the estimator as fussed nearly-isotonic approximator and denote it by $\hat{\bm{\beta}}^{FNI}(\bm{y}, \lambda_{F}, \lambda_{NI})$.

\section{The solution to the fused lasso nearly-isotonic signal approximator}
First, we consider fused nearly-isotonic regression, i.e. in (\ref{FLNIGim}) we assume that $\lambda_{L} = 0$. 

\begin{theorem}\label{sol_fnia}
For a fixed data vector $\bm{y} \in \mathbb{R}^{n}$ indexed by the index set $\mathcal{I}$ with the partial order relation $\preceq$ defined on $\mathcal{I}$ the solution to the fused nearly-isotonic problem in (\ref{FLNIGim}) is given by
\begin{equation}\label{FNI_sol}
\hat{\bm{\beta}}^{FNI}(\bm{y}, \lambda_{F}, \lambda_{NI}) = \bm{y} - D^{T} \hat{\bm{\nu}}(\lambda_{F}, \lambda_{NI})
\end{equation}
with
\begin{equation}\label{u_sol}
\hat{\bm{\nu}}(\bm{y}, \lambda_{F}, \lambda_{NI}) = \underset{\bm{\nu} \in \mathbb{R}^{n-1}}{\arg \min} \, \frac{1}{2}||\bm{y} - D^{T}\bm{\nu}||_{2}^{2} \quad \text{subject to} \quad - \lambda_{F} \bm{1} \leq \bm{\nu} \leq (\lambda_{F} + \lambda_{NI})\bm{1},
\end{equation}
\end{theorem}
where $D$ is the incidence matrix of the directed graph $G = (V, E)$, corresponding to $\preceq$ on  $\mathcal{I}$, $\bm{1}\in\mathbb{R}^{n}$ is the vector whose all elements are equal to $1$ and the notation $\bm{a} \leq \bm{b}$ for vectors $\bm{a},\bm{b} \in \mathbb{R}^{n}$ means $a_{i} \leq b_{i}$ for all $i = 1, \dots, n$.

\textbf{Proof.}
First, following the derivations of $\ell_{1}$ trend filtering and generalised lasso in \cite{kim2009ell_1}  and \cite{tibshirani2011solution}, respectively, we can write the optimization problem in (\ref{FNI}) in the following way
\begin{equation*}\label{}
\underset{\bm{\beta}, \bm{z}}{\text{minimize}} \, \frac{1}{2} ||\bm{y} - \bm{\beta}||_{2}^{2} +  \lambda_{F} ||\bm{z}||_{1} + \lambda_{NI}||\bm{z}||_{+} \quad \text{subject to} \quad D \bm{\beta} = \bm{z} \in \mathbb{R}^{n-1}.
\end{equation*}

Further, the Lagrangian is given by
\begin{equation}\label{Larg }
L(\bm{\beta}, \bm{z}, \bm{\nu}) = \frac{1}{2} ||\bm{y} - \bm{\beta}||_{2}^{2} + \lambda_{F} ||\bm{z}||_{1} + \lambda_{NI}||\bm{z}||_{+} + \bm{\nu}^{T}(D \bm{\beta} - \bm{z}),
\end{equation}
where $\bm{\nu} \in \mathbb{R}^{n-1}$ is a dual variable.

Note that 
\begin{equation*}\label{}
\underset{\bm{z}}{\min} \big(\lambda_{F} ||\bm{z}||_{1} +  \lambda_{NI}||\bm{z}||_{+} - \bm{\nu}^{T}\bm{z}\big) =
  \begin{cases}
    0, & \quad \text{if} \quad -\lambda_{F} \bm{1} \leq \bm{\nu} \leq (\lambda_{F} + \lambda_{NI})\bm{1}, \\
    -\infty, & \quad \text{otherwise},
  \end{cases}
\end{equation*}
and
\begin{equation*}\label{}
\underset{\bm{\beta}}{\min} \big( \frac{1}{2} ||\bm{y} - \bm{\beta}||_{2}^{2} + \bm{\nu}^{T} D \bm{\beta}\big) = -\frac{1}{2}\bm{\nu}^{T}DD^{T}\bm{\nu} + \bm{y}^{T}D^{T}\bm{\nu} = - \frac{1}{2} ||\bm{y} - D^{T}\bm{\nu}||_{2}^{2} + \frac{1}{2}\bm{y}^{T}\bm{y}.
\end{equation*}

Next, the dual function is given by
\begin{equation*}\label{}
\begin{aligned}
g(\nu) ={}& \underset{\bm{\beta}, \bm{z}}{\min} \, L(\bm{\beta}, \bm{z}, \bm{\nu}) =\\
  &\begin{cases}
    -\frac{1}{2} ||\bm{y} - D^{T}\bm{\nu}||_{2}^{2}  + \frac{1}{2}\bm{y}^{T}\bm{y}, & \quad \text{if} \quad -\lambda_{F} \bm{1} \leq \bm{\nu} \leq (\lambda_{F} + \lambda_{NI})\bm{1}, \\
    -\infty, & \quad \text{otherwise},
  \end{cases}
\end{aligned}
\end{equation*}
and, therefore, the dual problem is
\begin{equation*}\label{}
\hat{\bm{\nu}}(\bm{y}, \lambda_{F}, \lambda_{NI}) = \underset{\bm{\nu}}{\arg \max} \, g(\nu) \quad \text{subject to} \quad  -\lambda_{F} \bm{1} \leq \bm{\nu} \leq (\lambda_{F} + \lambda_{NI})\bm{1},
\end{equation*}
which is equivalent to
\begin{equation*}\label{}
\hat{\bm{\nu}}(\bm{y}, \lambda_{F}, \lambda_{NI}) = \underset{\bm{\nu}}{\arg \min} \, \frac{1}{2}||\bm{y} - D^{T}\bm{\nu}||_{2}^{2} \quad \text{subject to} \quad  -\lambda_{F} \bm{1} \leq \bm{\nu} \leq (\lambda_{F} + \lambda_{NI})\bm{1}.
\end{equation*}

Lastly, taking first derivative of Lagrangian $L(\bm{\beta}, \bm{z}, \bm{\nu})$ with respect to $\bm{\beta}$ we get the following relation between $\hat{\bm{\beta}}^{FNI}(\lambda_{F}, \lambda_{NI})$ and $\hat{\bm{\nu}}(\bm{y}, \lambda_{F}, \lambda_{NI})$
\begin{equation*}\label{}
\hat{\bm{\beta}}^{FNI}(\bm{y}, \lambda_{F}, \lambda_{NI}) = \bm{y} - D^{T} \hat{\bm{\nu}}(\bm{y},  \lambda_{F}, \lambda_{NI}).
\end{equation*}
\eop

Next, we provide the solution to the fused lasso nearly-isotonic regression.
\begin{theorem}\label{sol_FLNIa}
For a given vector $\bm{y}$ indexed by $\mathcal{I}$ the solution to the fused lasso nearly-isotonic signal approximator $\hat{\bm{\beta}}^{FLNI}(\bm{y},\lambda_{F}, \lambda_{L},\lambda_{NI})$ is given by soft thresholding the fused nearly-isotonic regression $\hat{\bm{\beta}}^{FNI}(\bm{y},\lambda_{F}, \lambda_{NI})$, i.e. 
\begin{equation}\label{FLNI_sol}
\hat{\beta}^{FLNI}_{\bm{i}}(\bm{y},\lambda_{F}, \lambda_{L},\lambda_{NI}) =   \begin{cases}
   \hat{\beta}^{FNI}_{\bm{i}}(\bm{y},\lambda_{F}, \lambda_{NI}) - \lambda_{L}, & \quad \text{if} \quad \hat{\beta}^{FNI}_{\bm{i}}(\bm{y},\lambda_{F}, \lambda_{NI}) \geq \lambda_{L}, \\
    0, & \quad \text{if} \quad |\hat{\beta}^{FNI}_{\bm{i}}(\bm{y},\lambda_{F}, \lambda_{NI})| \leq \lambda_{L},\\
    \hat{\beta}^{FNI}_{\bm{i}}(\bm{y},\lambda_{F}, \lambda_{NI}) + \lambda_{L}, & \quad \text{if} \quad \hat{\beta}^{FNI}_{\bm{i}}(\bm{y},\lambda_{F}, \lambda_{NI}) \leq -\lambda_{L},
  \end{cases}
\end{equation}
for ${\bm{i}}\in\mathcal{I}$.
\end{theorem}
\textbf{Proof.}
The proof is similar to the derivation of solution of the fused lasso in \cite{friedman2007pathwise}. Nevertheless, for  compliteness of the paper we provide the proof for $\hat{\bm{\beta}}^{FLNI}(\bm{y},\lambda_{F}, \lambda_{L},\lambda_{NI})$.

The subgradient equations (which are necessary and sufficient conditions for the solution of (\ref{FLNI})) for $\beta_{\bm{i}}$, with $\bm{i}\in\mathcal{I}$, are 
\begin{equation}\label{sgeq}
\begin{aligned}
g_{\bm{i}}(\lambda_{L}) ={}& -(y_{\bm{i}} - \beta_{\bm{i}}) + \lambda_{NI}(\sum_{\bm{j}:(\bm{i}, \bm{j})\in E} q_{\bm{i} , \bm{j}} - \sum_{\bm{j}:(\bm{j}, \bm{i})\in E}q_{\bm{j} , \bm{i}}) +\\
&\lambda_{F}(\sum_{\bm{j}:(\bm{i}, \bm{j})\in E}t_{\bm{i} , \bm{j}} - \sum_{\bm{j}:(\bm{j}, \bm{i})\in E}t_{\bm{j} , \bm{i}}) + \lambda_{L}s_{i} = 0,
\end{aligned}
\end{equation}
where 
\begin{equation}\label{qiti}
  q_{\bm{i}, \bm{j}} : \begin{cases}
     = 1, & \text{if } \, \beta_{\bm{i}} - \beta_{\bm{j}} > 0,\\
     = 0, & \text{if } \, \beta_{\bm{i}} - \beta_{\bm{j}} < 0,\\
    \in [0,1], & \text{if } \, \beta_{\bm{i}} =\beta_{\bm{j}},
  \end{cases} \quad \quad
  t_{\bm{i}, \bm{j}} : \begin{cases}
     = 1, & \text{if } \, \beta_{\bm{i}} - \beta_{\bm{j}} > 0,\\
     = -1, & \text{if } \, \beta_{\bm{i}} - \beta_{\bm{j}} < 0,\\
    \in [-1, 1], & \text{if } \, \beta_{\bm{i}} =\beta_{\bm{j}},
  \end{cases} 
\end{equation}

\begin{equation*}\label{}
s_{\bm{i}} : \begin{cases}
     = 1, & \text{if } \, \beta_{\bm{i}} > 0, \\
     = -1, & \text{if } \, \beta_{\bm{i}} < 0,\\
    \in [-1, 1], & \text{if } \, \beta_{\bm{i}} = 0.
  \end{cases}
\end{equation*}
Next, let $q_{\bm{i}, \bm{j}}(\lambda_{L})$, $t_{\bm{i}, \bm{j}}(\lambda_{L})$ and $s_{\bm{i}}(\lambda_{L})$ denote the values of the parameters defined above at some value of $\lambda_{L}$. Note, the values of $\lambda_{NI}$ and $\lambda_{F}$ are fixed. Therefore, if $\hat{\beta}_{\bm{i}}^{FLNI}(\bm{y},\lambda_{F}, 0, \lambda_{NI}) \neq 0$ for $s_{\bm{i}}(0)$ we have 
\begin{equation*}\label{}
s_{\bm{i}}(0) = \begin{cases}
     1, & \text{if } \, \hat{\beta}_{\bm{i}}^{FLNI}(\bm{y},\lambda_{F}, 0, \lambda_{NI}) > 0, \\
     -1, & \text{if } \, \hat{\beta}_{\bm{i}}^{FLNI}(\bm{y},\lambda_{F}, 0, \lambda_{NI}) < 0,\\
     \end{cases}
\end{equation*}
and for $\hat{\beta}_{\bm{i}}^{FLNI}(\bm{y}, \lambda_{F}, 0, \lambda_{NI}) = 0$ we can set $s_{\bm{i}}(0) = 0$.

Next, let $\hat{\bm{\beta}}^{ST}(\lambda_{L})$ denote the soft thresholding of $\hat{\bm{\beta}}^{FLNI}(\bm{y}, \lambda_{F}, 0, \lambda_{NI})$, i.e.

\begin{equation*}\label{}
\hat{\beta}_{\bm{i}}^{ST}(\lambda_{L}) =   \begin{cases}
   \hat{\beta}^{FLNI}_{\bm{i}}(\bm{y},\lambda_{F}, 0, \lambda_{NI}) - \lambda_{L}, & \quad \text{if} \quad \hat{\beta}_{\bm{i}}^{FLNI}(\bm{y}, \lambda_{F}, 0, \lambda_{NI}) \geq \lambda_{L}, \\
    0, & \quad \text{if} \quad |\hat{\beta}_{\bm{i}}^{FLNI}(\bm{y}, \lambda_{F}, 0,\lambda_{NI})| \leq \lambda_{L},\\
    \hat{\beta}_{\bm{i}}^{FLNI}(\bm{y}, \lambda_{F}, 0, \lambda_{NI}) + \lambda_{L}, & \quad \text{if} \quad \hat{\beta}_{\bm{i}}^{FLNI}(\bm{y},\lambda_{F}, 0, \lambda_{NI}) \leq -\lambda_{L}.
  \end{cases}
\end{equation*}
The goal is to prove that $\hat{\bm{\beta}}^{ST}(\lambda_{L})$ provides the solution to (\ref{FLNIG}). 

Note, analogously to the proof for the fused lasso estimator in Lemma A.1 at \cite{friedman2007pathwise}, if either $\hat{\beta}_{\bm{i}}^{ST}(\lambda_{L}) \neq 0$ or $\hat{\beta}_{\bm{j}}^{ST}(\lambda_{L}) \neq 0$, and $\hat{\beta}_{\bm{i}}^{ST}(\lambda_{L}) < \hat{\beta}_{\bm{j}}^{ST}(\lambda_{L})$ or $\hat{\beta}_{\bm{i}}^{ST}(\lambda_{L}) > \hat{\beta}_{\bm{j}}^{ST}(\lambda_{L})$,
then we also have $\hat{\beta}_{\bm{i}}^{ST}(0) < \hat{\beta}_{\bm{j}}^{ST}(0)$ or $\hat{\beta}_{\bm{i}}^{ST}(\lambda_{L}) > \hat{\beta}_{\bm{j}}^{ST}(\lambda_{L})$, respectively. Therefore, soft thresholding of $\hat{\bm{\beta}}^{FLNI}(\bm{y}, \lambda_{F}, 0, \lambda_{NI})$ does not change the ordering of these pairs and we have $q_{\bm{i}, \bm{j}}(\lambda_{L}) = q_{\bm{i},\bm{j}}(0)$ and $t_{\bm{i}, \bm{j}}(\lambda_{L}) = t_{\bm{i},\bm{j}}(0)$. Next, if for some $(\bm{i} , \bm{j})\in E$ we have $\hat{\beta}_{\bm{i}}^{ST}(\lambda_{L}) = \hat{\beta}_{\bm{j}}^{ST}(\lambda_{L}) = 0$, then $q_{\bm{i}, \bm{j}} \in [0,1]$ and $t_{\bm{i}, \bm{j}} \in [-1, 1]$, and, again, we can set $t_{\bm{i}, \bm{j}}(\lambda_{L}) = t_{\bm{i}, \bm{j}}(0)$, and $q_{\bm{i}, \bm{j}}(\lambda_{L}) = q_{\bm{i}, \bm{j}}(0)$.

Now let us insert $\hat{\beta}_{\bm{i}}^{ST}(\lambda_{L})$ into subgradient equations (\ref{sgeq}) and show that we can find $s_{\bm{i}}(\lambda_{L}) \in [0,1]$, for all $\bm{i}\in \mathcal{I}$.

First, assume that for some $\bm{i}$ we have $\hat{\beta}_{\bm{i}}^{FLNI}(\bm{y}, \lambda_{F}, 0, \lambda_{NI}) \geq \lambda_{L}$. Then

\begin{equation*}\label{}
\begin{aligned}
g_{\bm{i}}(\lambda_{L}) &={} -(y_{\bm{i}} - \hat{\beta}_{\bm{i}}^{FLNI}(\bm{y}, \lambda_{F}, 0, \lambda_{NI})) - \lambda_{L} + \lambda_{NI}(\sum_{\bm{j}:(\bm{i}, \bm{j})\in E} q_{\bm{i} , \bm{j}}(\lambda_{L}) - \sum_{\bm{j}:(\bm{j}, \bm{i})\in E}q_{\bm{j} , \bm{i}}(\lambda_{L}))\\ 
&+\lambda_{F}(\sum_{\bm{j}:(\bm{i}, \bm{j})\in E}t_{\bm{i} , \bm{j}}(\lambda_{L}) - \sum_{\bm{j}:(\bm{j}, \bm{i})\in E}t_{\bm{j} , \bm{i}}(\lambda_{L})) + \lambda_{L}s_{\bm{i}}(\lambda_{L})\\
&={} -(y_{\bm{i}} - \hat{\beta}_{\bm{i}}^{FLNI}(\bm{y}, \lambda_{F}, 0, \lambda_{NI})) + \lambda_{NI}(\sum_{\bm{j}:(\bm{i}, \bm{j})\in E} q_{\bm{i} , \bm{j}}(0) - \sum_{\bm{j}:(\bm{j}, \bm{i})\in E}q_{\bm{j} , \bm{i}}(0)) +\\
&\lambda_{F}(\sum_{\bm{j}:(\bm{i}, \bm{j})\in E}t_{\bm{i} , \bm{j}}(0) - \sum_{\bm{j}:(\bm{j}, \bm{i})\in E}t_{\bm{j} , \bm{i}}(0)) + \lambda_{L}s_{\bm{i}}(\lambda_{L}) - \lambda_{L} = 0.
\end{aligned}
\end{equation*}
Note, that 
\begin{equation*}\label{}
\begin{aligned}
  &-(y_{i} - \hat{\beta}_{i}^{FLNI}(\bm{y},\lambda_{F}, 0, \lambda_{NI})) + \lambda_{NI}(\sum_{\bm{j}:(\bm{i}, \bm{j})\in E} q_{\bm{i} , \bm{j}}(0) - \sum_{\bm{j}:(\bm{j}, \bm{i})\in E}q_{\bm{j} , \bm{i}}(0)) +\\
  &\lambda_{F}(\sum_{\bm{j}:(\bm{i}, \bm{j})\in E}t_{\bm{i} , \bm{j}}(0) - \sum_{\bm{j}:(\bm{j}, \bm{i})\in E}t_{\bm{j} , \bm{i}}(0)) = 0,
\end{aligned}
\end{equation*}
because $\hat{\bm{\beta}}^{FNI}(\bm{y}, \lambda_{F}, \lambda_{NI}) \equiv \hat{\bm{\beta}}^{FLNI}(\bm{y},\lambda_{F}, 0, \lambda_{NI})$. Therefore, if $s_{i}(\lambda_{L}) = \sign{\hat{\beta}_{\bm{i}}^{ST}(\lambda_{L})} = 1$, then $g_{\bm{i}}(\lambda_{L}) = 0$.  

The proof for the case when $\hat{\beta}_{\bm{i}}^{FLNI}(\bm{y}, \lambda_{F}, 0, \lambda_{NI}) \leq -\lambda_{L}$ is similar and one can show that $g_{\bm{i}}(\lambda_{L}) = 0$ if $s_{\bm{i}}(\lambda_{L}) = \sign{\hat{\beta}_{\bm{i}}^{ST}(\lambda_{L})} = -1$.

Second, assume that $|\hat{\beta}_{\bm{i}}^{FLNI}(\bm{y}, \lambda_{F}, 0, \lambda_{NI})| < \lambda_{L}$. Then, $\hat{\beta}_{\bm{i}}^{ST}(\lambda_{L}) = 0$, and

\begin{equation*}\label{}
\begin{aligned}
g_{\bm{i}}(\lambda_{L}) &={} -y_{\bm{i}} + \lambda_{NI}(\sum_{\bm{j}:(\bm{i}, \bm{j})\in E} q_{\bm{i} , \bm{j}}(\lambda_{L}) - \sum_{\bm{j}:(\bm{j}, \bm{i})\in E}q_{\bm{j} , \bm{i}}(\lambda_{L}))\\ 
&+\lambda_{F}(\sum_{\bm{j}:(\bm{i}, \bm{j})\in E}t_{\bm{i} , \bm{j}}(\lambda_{L}) - \sum_{\bm{j}:(\bm{j}, \bm{i})\in E}t_{\bm{j} , \bm{i}}(\lambda_{L})) + \lambda_{L}s_{\bm{i}}(\lambda_{L})\\
&={} -y_{\bm{i}} + \lambda_{NI}(\sum_{\bm{j}:(\bm{i}, \bm{j})\in E} q_{\bm{i} , \bm{j}}(0) - \sum_{\bm{j}:(\bm{j}, \bm{i})\in E}q_{\bm{j} , \bm{i}}(0)) +\\
&\lambda_{F}(\sum_{\bm{j}:(\bm{i}, \bm{j})\in E}t_{\bm{i} , \bm{j}}(0) - \sum_{\bm{j}:(\bm{j}, \bm{i})\in E}t_{\bm{j} , \bm{i}}(0)) + \lambda_{L}s_{\bm{i}}(\lambda_{L})  = 0.
\end{aligned}
\end{equation*}

Next, if we let $s_{\bm{i}}(\lambda_{L}) = \hat{\beta}_{\bm{i}}^{FLNI}(\bm{y},\lambda_{F}, 0, \lambda_{NI})/\lambda_{L}$, then, again, we have
\begin{equation*}\label{}
\begin{aligned}
  g_{\bm{i}}(\lambda_{L})&={}-(y_{i} - \hat{\beta}_{i}^{FLNI}(\bm{y},\lambda_{F}, 0, \lambda_{NI})) + \lambda_{NI}(\sum_{\bm{j}:(\bm{i}, \bm{j})\in E} q_{\bm{i} , \bm{j}}(0) - \sum_{\bm{j}:(\bm{j}, \bm{i})\in E}q_{\bm{j} , \bm{i}}(0)) +\\
  &\lambda_{F}(\sum_{\bm{j}:(\bm{i}, \bm{j})\in E}t_{\bm{i} , \bm{j}}(0) - \sum_{\bm{j}:(\bm{j}, \bm{i})\in E}t_{\bm{j} , \bm{i}}(0)) = 0,
\end{aligned}
\end{equation*}
Therefore, we have proved that $\hat{\bm{\beta}}^{FLNI}(\bm{y},\lambda_{F}, \lambda_{L},\lambda_{NI}) = \hat{\bm{\beta}}^{ST}(\lambda_{L})$.
\eop

\section{Properties of the fused lasso nearly-isotonic signal approximator}
We start with a proposition which shows how the solutions to the optimization problems (\ref{FG}), (\ref{NIG}) and (\ref{FLNIGim})  are related to each other. This result will be used in the next section to derive degrees of freedom of the fused lasso nearly-isotonic signal approximator.
\begin{proposition}\label{sol_rel}
For a fixed data vector $\bm{y}$ indexed by $\mathcal{I}$ and penalisation parameters $\lambda_{NI}$ and $\lambda_{F}$ the following relations between estimators $\hat{\bm{\beta}}^{F}$, $\hat{\bm{\beta}}^{NI}$ and $\hat{\bm{\beta}}^{FNI}$ hold
\begin{equation}\label{FNandI_sol}
\hat{\bm{\beta}}^{NI}(\bm{y}, \lambda_{NI}) = \hat{\bm{\beta}}^{F}(\bm{y} - \frac{\lambda_{NI}}{2} D^{T}\bm{1}, \frac{1}{2}\lambda_{NI}),
\end{equation}
\begin{equation}\label{FNI_F_NI_sol}
\hat{\bm{\beta}}^{FNI}(\bm{y}, \lambda_{F}, \lambda_{NI}) = \hat{\bm{\beta}}^{NI}(\bm{y} + \lambda_{F}D^{T}\bm{1}, \lambda_{NI} + 2\lambda_{F}) = \hat{\bm{\beta}}^{F}(\bm{y} - \frac{\lambda_{NI}}{2} D^{T}\bm{1}, \frac{1}{2}\lambda_{NI} + \lambda_{F})
\end{equation}
and
\begin{equation}\label{FNI_FL_NI_sol}
\hat{\bm{\beta}}^{FLNI}(\bm{y}, \lambda_{F}, \lambda_{L}, \lambda_{NI}) = \hat{\bm{\beta}}^{FL}(\bm{y}- \frac{\lambda_{NI}}{2} D^{T}\bm{1}, \frac{1}{2}\lambda_{NI} + \lambda_{F}, \lambda_{L}),
\end{equation}
where $D$ is the oriented incidence matrix for the graph $G = (V, E)$ corresponding to the partial order relation $\preceq$ on $\mathcal{I}$. 
\end{proposition}

\textbf{Proof.}
First, from \cite{tibshirani2011nearly} the solution to the nearly-isotonic problem is given by 
\begin{equation*}\label{}
\hat{\bm{\beta}}^{NI}(\bm{y}, \lambda_{NI}) = \bm{y} - D^{T} \hat{\bm{v}}(\bm{y}, \lambda_{NI}),
\end{equation*}
with
\begin{equation*}\label{}
\hat{\bm{v}}(\bm{y}, \lambda_{NI}) = \underset{\bm{v} \in \mathbb{R}^{n-1}}{\arg \min} \, \frac{1}{2}||\bm{y} - D^{T}\bm{v}||_{2}^{2} \quad \text{subject to} \quad   \bm{0} \leq \bm{v} \leq \lambda_{NI}\bm{1},
\end{equation*}
and from \cite{tibshirani2011solution} it follows
\begin{equation*}\label{} 
\hat{\bm{\beta}}^{F}(\bm{y}, \lambda_{F}) = \bm{y} - D^{T} \hat{\bm{w}}(\bm{y}, \lambda_{F}),
\end{equation*}
with
\begin{equation*}\label{}
\hat{\bm{w}}(\bm{y}, \lambda_{F}) = \underset{\bm{w} \in \mathbb{R}^{n-1}}{\arg \min} \, \frac{1}{2}||\bm{y} - D^{T}\bm{w}||_{2}^{2} \quad \text{subject to} \quad  -\lambda_{F} \bm{1} \leq \bm{w} \leq \lambda_{F}\bm{1}.
\end{equation*}

Second, let us introduce a new variable $\bm{v}^{*} = \bm{v} - \frac{\lambda_{NI}}{2} \bm{1}$. Then
\begin{equation*}\label{}
\hat{\bm{\beta}}^{NI}(\bm{y}, \lambda_{NI}) = \bm{y} -D^{T}\frac{\lambda_{NI}}{2} \bm{1} - D^{T} \hat{\bm{v}}^{*}(\bm{y}, \lambda_{NI}),
\end{equation*}
where
\begin{equation*}\label{}
\hat{\bm{v}}^{*}(\bm{y}, \lambda_{NI}) = \underset{\bm{v}^{*} \in \mathbb{R}^{n-1}}{\arg \min} \, \frac{1}{2}||\bm{y} - D^{T}\frac{\lambda_{NI}}{2} \bm{1} - D^{T}\bm{v}^{*}||_{2}^{2} \quad \text{subject to} \quad   -\frac{\lambda_{NI}}{2} \bm{1} \leq \bm{v}^{*} \leq \frac{\lambda_{NI}}{2} \bm{1}.
\end{equation*}
Therefore, we have proved that $\hat{\bm{\beta}}^{NI}(\bm{y}, \lambda_{NI}) = \hat{\bm{\beta}}^{F}(\bm{y} - \frac{\lambda_{NI}}{2} D^{T}\bm{1}, \frac{1}{2}\lambda_{NI})$. 

The proof for the fused lasso nearly-isotonic estimator is the same with the change of variable $\bm{u}^{*} = \bm{u} + D^{T}\lambda_{F} \bm{1}$ in (\ref{FNI_sol}) and (\ref{u_sol}) for the proof of the first equality in (\ref{FNI_F_NI_sol}) and with $\bm{u}^{*} = \bm{u} - \frac{\lambda_{NI}}{2} \bm{1}$ for the second equality.

Next, we prove the result for the case of fused lasso nearly-isotonic approximator. From Theorem \ref{sol_FLNIa} we have
\begin{equation*}\label{}
\hat{\beta}^{FLNI}_{\bm{i}}(\bm{y},\lambda_{F}, \lambda_{L},\lambda_{NI}) =   \begin{cases}
   \hat{\beta}^{FNI}_{\bm{i}}(\bm{y},\lambda_{F}, \lambda_{NI}) - \lambda_{L}, & \quad \text{if} \quad \hat{\beta}^{FNI}_{\bm{i}}(\bm{y},\lambda_{F}, \lambda_{NI}) \geq \lambda_{L}, \\
    0, & \quad \text{if} \quad |\hat{\beta}^{FNI}_{\bm{i}}(\bm{y},\lambda_{F}, \lambda_{NI})| \leq \lambda_{L},\\
    \hat{\beta}^{FNI}_{\bm{i}}(\bm{y},\lambda_{F}, \lambda_{NI}) + \lambda_{L}, & \quad \text{if} \quad \hat{\beta}^{FNI}_{\bm{i}}(\bm{y},\lambda_{F}, \lambda_{NI}) \leq -\lambda_{L},
  \end{cases}
\end{equation*}
for $\bm{i} \in \mathcal{I}$.

Further, using (\ref{FNI_F_NI_sol}) we have
\begin{equation*}\label{}
\hat{\beta}^{FLNI}_{\bm{i}}(\bm{y},\lambda_{F}, \lambda_{L},\lambda_{NI}) = 
   \hat{\beta}^{F}_{\bm{i}}(\bm{y} - \frac{\lambda_{NI}}{2} D^{T}\bm{1}, \frac{1}{2}\lambda_{NI} + \lambda_{F}) - \lambda_{L},
\end{equation*}
if $\hat{\beta}^{F}_{\bm{i}}(\bm{y} - \frac{\lambda_{NI}}{2} D^{T}\bm{1}, \frac{1}{2}\lambda_{NI} + \lambda_{F}) \geq \lambda_{L}$,
\begin{equation*}\label{}
\hat{\beta}^{FLNI}_{\bm{i}}(\bm{y},\lambda_{F}, \lambda_{L},\lambda_{NI}) = 0,
\end{equation*}
if $|\hat{\beta}^{F}_{\bm{i}}(\bm{y} - \frac{\lambda_{NI}}{2} D^{T}\bm{1}, \frac{1}{2}\lambda_{NI} + \lambda_{F})| \leq \lambda_{L}$,
\begin{equation*}\label{}
\hat{\beta}^{FLNI}_{\bm{i}}(\bm{y},\lambda_{F}, \lambda_{L},\lambda_{NI}) = 
   \hat{\beta}^{F}_{\bm{i}}(\bm{y} - \frac{\lambda_{NI}}{2} D^{T}\bm{1}, \frac{1}{2}\lambda_{NI} + \lambda_{F}) + \lambda_{L},
\end{equation*}
if $\hat{\beta}^{F}_{\bm{i}}(\bm{y} - \frac{\lambda_{NI}}{2} D^{T}\bm{1}, \frac{1}{2}\lambda_{NI} + \lambda_{F}) \leq -\lambda_{L}$.

Therefore, we obtain 
\begin{equation*}\label{}
\begin{aligned}
\hat{\bm{\beta}}^{FLNI}(\bm{y}, \lambda_{F}, \lambda_{L}, \lambda_{NI}) &={}\\
\underset{\bm{\beta} \in \mathbb{R}^{n}}{\arg \min} \, \frac{1}{2} ||\bm{y} &- \frac{\lambda_{NI}}{2} D^{T}\bm{1} - \bm{\beta}||_{2}^{2} + (\frac{1}{2}\lambda_{NI} + \lambda_{F}) ||D\bm{\beta}||_{1} + \lambda_{L} ||\bm{\beta}||_{1} \equiv\\
\hat{\bm{\beta}}^{FL}(\bm{y} &- \frac{\lambda_{NI}}{2} D^{T}\bm{1}, \frac{1}{2}\lambda_{NI} + \lambda_{F}, \lambda_{L}).
\end{aligned}
\end{equation*}
\eop

Next, we prove that, analogously to fussed lasso and nearly-isotonic regression, as one of the penalization parameters increases the constant regions in the solution $\hat{\bm{\beta}}^{FLNI}$ can only be joined together and not split apart. We prove this result only for one dimensional monotonic order, and the general case is an open question. This result could be potentially useful in the future research for the solution path of fussed lasso nearly-isotonic approximator.
\begin{proposition}
Let $\mathcal{I} = \{1, \dots, n\}$ with a natural order defined on it. Next, let $\bm{\lambda} = (\lambda_{F}, \lambda_{L}, \lambda_{NI})$ and $\bm{\lambda}^{*}= (\lambda_{F}^{*}, \lambda_{L}^{*}, \lambda_{NI}^{*})$ are the triples of penalisation parameters such that one of the elements of $\bm{\lambda}^{*}$ is greater than the corresponding element in $\bm{\lambda}$, while two others are the same. Next, assume that for some $i$ the solution $\hat{\bm{\beta}}^{FLNI}(\bm{y},\bm{\lambda})$ satisfies \begin{equation*}\label{}
\hat{\beta}_{i}^{FLNI}(\bm{y},\bm{\lambda}) = \hat{\beta}_{i+1}^{FLNI}(\bm{y},\bm{\lambda}).
\end{equation*}
Then for $\bm{\lambda}^{*}$ we have
\begin{equation*}\label{}
\hat{\beta}_{i}^{FLNI}(\bm{y},\bm{\lambda}^{*}) = \hat{\beta}_{i+1}^{FLNI}(\bm{y},\bm{\lambda}^{*}).
\end{equation*}
\end{proposition}
\textbf{Proof.}

\textbf{Case 1: $\lambda_{NI}$ and $\lambda_{F}$ are fixed and $\lambda_{L}^{*} > \lambda_{L}$.}
The result of the proposition for this case follows directly from Theorem \ref{sol_FLNIa}.

\textbf{Case 2: $\lambda_{F}$ and $\lambda_{L}$ are fixed and $\lambda_{NI}^{*} > \lambda_{NI}$.}
Let us consider the fused nearly-isotonic regression and write the subgradient equations
\begin{equation*}\label{}
g_{i}(\lambda_{NI}) = -(y_{i} - \beta_{i}) + \lambda_{NI}(q_{i}(\lambda_{NI}) - q_{i-1}(\lambda_{NI})) + \lambda_{F}(t_{i}(\lambda_{NI}) - t_{i-1}(\lambda_{NI})) = 0,
\end{equation*}
where $q_{i}$ and $t_{i}$, with $i = 1, \dots, n$, are defined in (\ref{qiti}), and, analogously to the proof of Theorem \ref{sol_FLNIa}, $q(\lambda_{NI})$, $t(\lambda_{NI})$ denote the values of the parameters defined above at some value of $\lambda_{NI}$.

Assume that for $\lambda_{NI}$ in the solution $\hat{\bm{\beta}}^{FNI}(\bm{y}, \lambda_{F}, \lambda_{NI})$ we have a following constant region 
\begin{equation}\label{const-reg-FNI}
\hat{\beta}_{j-1}^{FNI}(\bm{y}, \lambda_{F}, \lambda_{NI}) \neq \hat{\beta}_{j}^{FNI}(\bm{y}, \lambda_{F}, \lambda_{NI}) = \dots = \hat{\beta}_{j+k}^{FNI}(\bm{y},\lambda_{F}, \lambda_{NI}) \neq \hat{\beta}_{j+k +1}^{FNI}(\bm{y}, \lambda_{F}, \lambda_{NI}),
\end{equation}
and in the similar way as in \cite{tibshirani2011nearly} for $\lambda_{NI}^{*}$ we consider the subset of the subgradient equations
\begin{equation}\label{sgesbs}
g_{i}(\lambda_{NI}) = -(y_{i} - \beta_{i}) + \lambda_{NI}^{*}(q_{i}(\lambda_{NI}^{*}) - q_{i-1}(\lambda_{NI}^{*})) + \lambda_{F}(t_{i}(\lambda_{NI}^{*}) - t_{i-1}(\lambda_{NI}^{*})) = 0,
\end{equation}
with $i = j, \dots, k$, and show that there exists the solution for which (\ref{const-reg-FNI}) holds, $q_{i} \in [0, 1]$ and $t_{i} \in [-1, 1]$.

Note first that as $\lambda_{NI}$ increases, (\ref{const-reg-FNI}) holds until the merge with other groups happens, which means that $q_{j-1}, q_{j+k} \in \{0, 1\}$ and $t_{j-1}, t_{j+k} \in \{-1, 1\}$ will not change their values until the merge of this constant region. Also, as it follows from (\ref{qiti}), for $i \in [j, j+ k]$ the value of $t_{i}$ is in $[-1 ,1]$. Therefore, without any violation of the restrictions on $t_{i}$ we can assume that $t_{i}(\lambda_{NI}^{*}) = t_{i}(\lambda)$ for any $i \in [j, j + k - 1]$.

Next, taking pairwise differences between subgradient equations for $\lambda_{NI}$ we have
\begin{equation*}\label{}
\lambda_{NI}A\tilde{\bm{q}}(\lambda_{NI}) + \lambda_{F}A\tilde{\bm{t}}(\lambda_{NI}) = D\tilde{\bm{y}} + \lambda_{NI}\bm{c}(\lambda_{NI}) + \lambda_{F}\bm{d}(\lambda_{NI}),
\end{equation*}
where $D$ is displayed at Figure \ref{mgrmtr}, 
\begin{equation}\label{Amatrix}
A = \begin{bmatrix}
2 & -1 & 0 & \dots & 0 & 0 & 0 \\
-1 & 2 & -1 & \dots & 0 & 0 & 0 \\
\vdots & \vdots  & \vdots  & \vdots & \vdots  & \vdots & \vdots \\
0 & 0 & 0 & \dots & -1 & 2 & -1\\
0 & 0 & 0 & \dots & 0 & -1 & 2 
\end{bmatrix},
\end{equation}
and $\tilde{\bm{q}}(\lambda_{NI}) = (q_{j}(\lambda_{NI}), \dots, q_{j + k - 1}(\lambda_{NI})) $, $\tilde{\bm{t}}(\lambda_{NI}) = (t_{j}(\lambda_{NI}), \dots, t_{j + k - 1}(\lambda_{NI})) $, $\tilde{\bm{y}} = (y_{j}, \dots, y_{j + k})$, $\bm{c}(\lambda_{NI}) = (q_{j-1}(\lambda_{NI}), 0, \dots, 0, q_{j+k}(\lambda_{NI}))$, and $\bm{d}(\lambda_{NI}) = (t_{j-1}(\lambda_{NI}), 0, \dots, 0,  t_{j+k}(\lambda_{NI}))$. 

Since $A$ is invertible we have
\begin{equation*}\label{}
\lambda_{NI}\tilde{\bm{q}}(\lambda_{NI}) + \lambda_{F}\tilde{\bm{t}}(\lambda_{NI}) = A^{-1}D\tilde{\bm{y}} + \lambda_{NI}A^{-1}\bm{c}(\lambda_{NI}) + \lambda_{F}A^{-1}\bm{d}(\lambda_{NI}),
\end{equation*}
and, since $\tilde{\bm{q}}(\lambda_{NI})$ and $\tilde{\bm{t}}(\lambda_{NI})$ provide the solution to the subgradient equations (\ref{sgesbs}), then
\begin{equation}\label{l1qt}
-\lambda_{F} \leq \lambda_{NI}\tilde{\bm{q}}(\lambda_{NI}) + \lambda_{F}\tilde{\bm{t}}(\lambda_{NI}) \leq \lambda_{NI} + \lambda_{F}.
\end{equation}

Next, as pointed out at \cite{friedman2007pathwise} and \cite{tibshirani2011nearly} 
\begin{equation*}\label{}
(A^{-1})_{i,1} = (n-i+1)/(n+1) \quad \text{and} \quad (A^{-1})_{i,n} = i/(n+1),
\end{equation*}
then, one can show that
\begin{equation}\label{l1cd}
-\lambda_{F}\bm{1} \preceq \lambda_{NI}A^{-1}\bm{c}(\lambda_{NI}) + \lambda_{F}A^{-1}\bm{d}(\lambda_{NI}) \preceq \lambda_{NI}\bm{1} + \lambda_{F}\bm{1}.
\end{equation}

Further, let us consider the case of $\lambda_{NI}^{*} > \lambda_{NI}$. Then we have
\begin{equation*}\label{}
\lambda_{NI}^{*}\tilde{\bm{q}}(\lambda_{NI}^{*}) + \lambda_{F}\tilde{\bm{t}}(\lambda_{NI}^{*}) = A^{-1}D\tilde{\bm{y}} + \lambda_{NI}^{*}A^{-1}\bm{c}(\lambda_{NI}^{*}) + \lambda_{F}A^{-1}\bm{d}(\lambda_{NI}^{*}).
\end{equation*}
Recall, above we set $\tilde{\bm{t}}(\lambda_{NI}^{*}) = \tilde{\bm{t}}(\lambda_{NI})$, and $q_{j-1}, q_{j+k}, t_{j-1}$ and $t_{j+k}$ does not change their values until the merge, which means that $\bm{c}(\lambda_{NI}^{*}) = \bm{c}(\lambda_{NI})$, and $\bm{d}(\lambda_{NI}^{*}) = \bm{d}(\lambda_{NI})$. 

Therefore, the subgradient equations for $\lambda_{NI}^{*}$ can be written as 
\begin{equation*}\label{}
\lambda_{NI}^{*}\tilde{\bm{q}}(\lambda_{NI}^{*}) + \lambda_{F}\tilde{\bm{t}}(\lambda_{NI}) = A^{-1}D\tilde{\bm{y}} + \lambda_{NI}^{*}A^{-1}\bm{c}(\lambda_{NI}) + \lambda_{F}A^{-1}\bm{d}(\lambda_{NI}).
\end{equation*}
Next, since the term $A^{-1}D\tilde{\bm{y}}$ is not changed, $-\lambda_{F} \leq \lambda_{F}\tilde{\bm{t}}(\lambda_{NI}) \leq \lambda_{F}$, and 
\begin{equation*}\label{}
-\lambda_{F}\bm{1} \preceq \lambda_{NI}^{*}A^{-1}\bm{c}(\lambda_{NI}) + \lambda_{F}A^{-1}\bm{d}(\lambda_{NI}) \preceq \lambda_{NI}^{*}\bm{1} + \lambda_{F}\bm{1},
\end{equation*}
then we have 
\begin{equation*}\label{}
\bm{0} \preceq \tilde{\bm{q}}(\lambda_{NI}^{*}) \preceq \bm{1}.
\end{equation*}
Therefore we proved that $\hat{\beta}_{i}^{FNI}(\bm{y},\bm{\lambda}^{*}) = \hat{\beta}_{i+1}^{FNI}(\bm{y},\bm{\lambda}^{*})$. Since $\hat{\beta}_{i}^{FLNI}(\bm{y},\bm{\lambda}^{*})$ is given by soft thresholding of $\hat{\beta}_{i}^{FNI}(\bm{y},\bm{\lambda}^{*})$, then $\hat{\beta}_{i}^{FLNI}(\bm{y},\bm{\lambda}^{*}) = \hat{\beta}_{i+1}^{FLNI}(\bm{y},\bm{\lambda}^{*})$ for $i \in [j, k]$.

\textbf{Case 3: $\lambda_{NI}$ and $\lambda_{L}$ are fixed and $\lambda_{F}^{*} > \lambda_{F}$.} The proof for this case is virtually identical to the proof for the Case 2. In this case we assume that $q_{i}(\lambda_{F}^{*}) = q_{i}(\lambda_{2 })$ for any $i \in [j, j + k-1]$. Next, $q_{j-1}, q_{j+k}, t_{j-1}$ and $t_{j+k}$ do not change their values until the merge, which, again, means that $\bm{c}(\lambda_{F}^{*}) = \bm{c}(\lambda_{F})$, and $\bm{d}(\lambda_{F}^{*}) = \bm{d}(\lambda_{F})$. Finally, we can show that 
\begin{equation*}\label{}
-\bm{1} \preceq \tilde{\bm{t}}(\lambda_{F}^{*}) \preceq \bm{1}.
\end{equation*}
\eop

\section{Degrees of freedom}
In this section we discuss the estimation of degrees of freedom for fused nearly-isotonic regression and fused lasso nearly isotonic signal approximator. Let us consider the following nonparametric model 
\begin{equation*}\label{}
\bm{Y} = \bm{\mathring{\beta}} + \bm{\varepsilon},
\end{equation*}
where $\bm{\mathring{\beta}}\in\mathbb{R}^{n}$ is the unknown signal, and the error term $\bm{\varepsilon}\in\mathcal{N}(\bm{0}, \sigma^{2}\bm{I})$, with $\sigma < \infty$. 

The degrees of freedom is a measure of complexity of the estimator, and following \citet{efron1986biased},  for the fixed values of $\lambda_{F}$, $\lambda_{L}$ and $\lambda_{Ni}$ the degrees of freedom of $\hat{\bm{\beta}}^{FNI}$ and $\hat{\bm{\beta}}^{FLNI}$ are given by
\begin{equation}\label{dfFNI}
df(\hat{\bm{\beta}}^{FNI}(\bm{Y}, \lambda_{F}, \lambda_{NI})) = \frac{1}{\sigma^{2}} \sum_{i=1}^{n}\Cov[\hat{\beta}^{FNI}_{i}(\bm{Y}, \lambda_{F}, \lambda_{NI}), Y_{i}]
\end{equation}
and
\begin{equation}\label{dfFLNI}
df(\hat{\bm{\beta}}^{FLNI}(\bm{Y}, \lambda_{F}, \lambda_{L}, \lambda_{NI})) = \frac{1}{\sigma^{2}} \sum_{i=1}^{n}\Cov[\hat{\beta}^{FLNI}_{i}(\bm{Y}, \lambda_{F},\lambda_{L}, \lambda_{NI}), Y_{i}].
\end{equation}

The next theorem provides the unbiased estimators of the degrees of freedom $df(\hat{\bm{\beta}}^{FNI})$ and $df(\hat{\bm{\beta}}^{FLNI})$.
\begin{theorem}\label{estdf}
For the fixed values of $\lambda_{F}$, $\lambda_{L}$ and $\lambda_{Ni}$ let 
\begin{equation*}\label{}
K^{FNI}(\bm{y}, \lambda_{F}, \lambda_{NI}) = \#\{\text{fused groups in } \hat{\bm{\beta}}^{FNI}(\bm{y}, \lambda_{F}, \lambda_{NI})\},
\end{equation*}
and 
\begin{equation*}\label{}
K^{FLNI}(\bm{y}, \lambda_{F}, \lambda_{L}, \lambda_{NI}) = \#\{\text{non-zero fused groups in } \hat{\bm{\beta}}^{FLNI}(\bm{y}, \lambda_{F}, \lambda_{L}, \lambda_{NI})\}.
\end{equation*}
Then we have 
\begin{equation*}\label{}
\mathbb{E}[K^{FNI}(\bm{Y}, \lambda_{F}, \lambda_{NI})] = df(\hat{\bm{\beta}}^{FNI}(\bm{Y}, \lambda_{F}, \lambda_{NI})), 
\end{equation*}
and
\begin{equation*}\label{}
\mathbb{E}[K^{FLNI}(\bm{Y}, \lambda_{F}, \lambda_{L}, \lambda_{NI})] = df(\hat{\bm{\beta}}^{FLNI}(\bm{Y}, \lambda_{F}, \lambda_{L}, \lambda_{NI})).
\end{equation*}
\end{theorem}
\textbf{Proof.}
First, for the fused estimator $\hat{\bm{\beta}}^{F}(\bm{y}, \lambda_{F})$ let 
\begin{equation*}\label{}
K^{F}(\bm{y}, \lambda_{F}) = \#\{\text{fused groups in } \hat{\bm{\beta}}^{F}(\bm{y}, \lambda_{F})\}.
\end{equation*}
Then, as it follows from \citet{tibshirani2011solution}, for $\hat{\bm{\beta}}^{F}(\bm{y}, \lambda_{F})$ we have 
\begin{equation*}\label{}
\mathbb{E}[K^{F}(\bm{Y}, \lambda_{F})] = df(\hat{\bm{\beta}}^{F}(\bm{Y}, \lambda_{F})).
\end{equation*}

Next, from Proposition \ref{sol_rel}, it follows 
\begin{equation*}\label{}
\hat{\bm{\beta}}^{FNI}(\bm{y}, \lambda_{F}, \lambda_{NI}) = \hat{\bm{\beta}}^{F}(\bm{y} - \frac{\lambda_{NI}}{2} D^{T}\bm{1}, \frac{1}{2}\lambda_{NI} + \lambda_{F}).
\end{equation*}
Therefore, using the property of covariance we have

\begin{equation*}\label{}
\begin{aligned}
df(\hat{\bm{\beta}}^{FNI}(&\bm{Y}, \lambda_{F}, \lambda_{NI}))  =\sum_{i=1}^{n}\Cov[\hat{\beta}^{FNI}_{i}(\bm{Y}, \lambda_{F}, \lambda_{NI}), Y_{i}] =\\ &\sum_{i=1}^{n}\Cov[\hat{\beta}^{F}_{i}(\bm{Y} - \frac{\lambda_{NI}}{2} D^{T}\bm{1}, \frac{1}{2}\lambda_{NI} + \lambda_{F}), Y_{i}] =\\ 
&\sum_{i=1}^{n}\Cov[\hat{\beta}^{F}_{i}(\bm{Y} - \frac{\lambda_{NI}}{2} D^{T}\bm{1}, \frac{1}{2}\lambda_{NI} + \lambda_{F}), Y_{i} - \frac{\lambda_{NI}}{2} [D^{T}\bm{1}]_{i}]=\\
&\mathbb{E}[K^{F}(\bm{Y} - \frac{\lambda_{NI}}{2} D^{T}\bm{1}, \frac{1}{2}\lambda_{NI} + \lambda_{F})] \equiv \mathbb{E}[K^{FNI}(\bm{Y}, \lambda_{F} , \lambda_{NI})],
\end{aligned}
\end{equation*}
where $[\bm{a}]_{i}$ denotes $i$-th element in the vector $\bm{a}\in\mathbb{R}^{n}$. 

Next, we prove the result for the fused lasso nearly-isotonic approximator. From Proposition \ref{sol_rel} we have 
\begin{equation*}\label{}
\hat{\bm{\beta}}^{FLNI}(\bm{y}, \lambda_{F}, \lambda_{L}, \lambda_{NI}) = \hat{\bm{\beta}}^{FL}(\bm{y}- \frac{\lambda_{NI}}{2} D^{T}\bm{1}, \frac{1}{2}\lambda_{NI} + \lambda_{F}, \lambda_{L}).
\end{equation*}
Next, for the fused lasso $\hat{\bm{\beta}}^{FL}(\bm{y}, \lambda_{F}, \lambda_{L})$ defined in (\ref{FL}) let
\begin{equation*}\label{}
K^{FL}(\bm{y}, \lambda_{F}, \lambda_{L}) = \#\{\text{non-zero fused groups in } \hat{\bm{\beta}}^{FL}(\bm{y}, \lambda_{F}, \lambda_{L})\},
\end{equation*}
and from \citet{tibshirani2011solution} it follows
\begin{equation*}\label{}
\mathbb{E}[K^{FL}(\bm{Y}, \lambda_{F}, \lambda_{L})] = df(\hat{\bm{\beta}}^{FL}(\bm{Y}, \lambda_{F}, \lambda_{L})).
\end{equation*}

Further, again, using the property of the covariance, we have
\begin{equation*}\label{}
\begin{aligned}
df(\hat{\bm{\beta}}^{FLNI}&(\bm{Y}, \lambda_{F}, \lambda_{L}, \lambda_{NI})) ={} \sum_{i=1}^{n}\Cov[\hat{\beta}^{FLNI}_{i}(\bm{Y}, \lambda_{F},\lambda_{L}, \lambda_{NI}), Y_{i}] =\\
&\sum_{i=1}^{n}\Cov[\hat{\beta}^{FL}_{i}(\bm{Y}- \frac{\lambda_{NI}}{2} D^{T}\bm{1}, \frac{1}{2}\lambda_{NI} + \lambda_{F}, \lambda_{L}), Y_{i}] =\\ 
&\sum_{i=1}^{n}\Cov[\hat{\beta}^{FL}_{i}(\bm{Y} - \frac{\lambda_{NI}}{2} D^{T}\bm{1}, \frac{1}{2}\lambda_{NI} + \lambda_{F},\lambda_{L}), Y_{i} - \frac{\lambda_{NI}}{2} [D^{T}\bm{1}]_{i}]=\\
&\mathbb{E}[K^{FL}(\bm{Y} - \frac{\lambda_{NI}}{2} D^{T}\bm{1}, \frac{1}{2}\lambda_{NI} + \lambda_{F},\lambda_{L})]\equiv \mathbb{E}[K^{FLNI}(\bm{Y}, \lambda_{F}, \lambda_{L}, \lambda_{NI})].
\end{aligned}
\end{equation*}

Lastly, we note that the proof for the unbiased estimator of the degrees of freedom for nearly-isotonic regression, given in \citet{tibshirani2011nearly}, can be done in the same way as in the current proof, using the relation (\ref{FNandI_sol}) and, again, the result of the paper \citet{tibshirani2011solution} for the fusion estimator $\hat{\bm{\beta}}^{FLNI}(\bm{Y}, \lambda_{F})$.
\eop

We can use the estimate of degrees of freedom for unbiased estimation of the true risk $\mathbb{E}[\sum_{i=1}^{n}(\mathring{\beta}_{i}- \hat{\beta}^{FLNI}_{i}(\bm{Y}, \lambda_{F}, \lambda_{L}, \lambda_{NI}))^2]$, which is given by $\hat{C}_{p}$ statistic
\begin{equation*}\label{}
\hat{C}_{p}(\lambda_{F}, \lambda_{L}, \lambda_{NI}) = \sum_{i=1}^{n}(y_{i} - \hat{\beta}^{FLNI}_{i}(\bm{y}, \lambda_{F}, \lambda_{L}, \lambda_{NI}))^2 - n\sigma^{2} + 2\sigma^{2}K^{FLNI}(\bm{Y}, \lambda_{F}, \lambda_{L}, \lambda_{NI}).
\end{equation*}

\section{Acknowledgments}
This work was partially supported by the Wallenberg AI, Autonomous Systems and Software Program (WASP) funded by the Knut and Alice Wallenberg Foundataion. 

\bigskip

\appendix
\section{Appendix}

\bibliographystyle{agsm}

\bibliography{Bibliography-MM-MC}
\end{document}